\documentclass[12pt]{amsart}
\usepackage{amsmath,amssymb}
\usepackage{color}

\theoremstyle{plain}
\newtheorem{thm}{Theorem}[section]
\newtheorem{lem}[thm]{Lemma}
\newtheorem{cor}[thm]{Corollary}
\newtheorem{prop}[thm]{Proposition}

\newtheorem{question}[thm]{Question}

\theoremstyle{definition}

\newtheorem{example}[thm]{Example}

\begin{document}

\title{\bf Simple polynomial equations over $(m \times m)$-matrices}

\author{Vitalij A.~Chatyrko and Alexandre Karassev}


\begin{abstract}
Let $m$ be any integer $\geq 3$. We consider the polynomial equation
$$X^n + a_{n-1}\cdot X^{n-1} + \dots + a_1 \cdot X + a_0 \cdot I = O,$$
over $(m \times m)$-matrices $X$ with the real entries, where $I$ is the identity matrix, $O$ is the null matrix, $a_i \in \mathbb R$
for each $i$ and $n \geq 1$. We discuss its solution set $S$ supplied with the natural Euclidean topology.
In particular, we describe the solution set $S$ for $m=3$ and calculate its dimension.

\end{abstract}



\maketitle

\section{Introduction}

Let $m$ be any integer $\geq 2$. As usual, by $M_m(\mathbb R)$  we denote the algebra of real $(m\times m)$ matrices, and by $GL_m (\mathbb R)$ the group of invertible $(m\times m)$  matrices. We equip $M_m (\mathbb R)$ with any of equivalent norms that make  $M_m(\mathbb R)$ a topological algebra, homeomorphic to $\mathbb R^{m^2}.$

Consider a polynomial  $$F(X) = X^n + A_{n-1}X^{n-1} + \dots + A_1 X + A_0$$ of degree  $n \geq 2,$
where $X$ and $A_i$, $i=0,1,\dots, n-1,$ are from  $M_m(\mathbb R)$, and the equation
$$F(X) = O \eqno{(1)},$$
where $O$ is the null $(m\times m)$-matrix.

Any matrix $B \in M_m(\mathbb R)$ such that $F(B) = O$ is called {\it a solution} of $(1)$. {\it The solution set} of $(1)$ consists of all solutions of $(1)$.

In \cite{W} Wilson demonstrated a method showing that the solution set of  $(1)$ with $n=2$
over $M_2(\mathbb R)$
can consist of exactly $0, 1, 2, 3, 4, 5, 6$ elements or it can be infinite (see also \cite{G} and \cite{S}).  The method uses eigenvalues and eigenvectors and  is due to Fuchs and Schwarz \cite{FS}. 

One can simplify the polynomial  $F(X)$ by setting $A_k = a_k \cdot I$ for each $0 \leq k \leq n-1$, where each $a_k$ is a real number. Denote the polynomial  $F(X)$ in this case by $F_s(X)$ and note that $$F_s(X) = X^n + a_{n-1}\cdot X^{n-1} + \dots + a_1 \cdot X + a_0 \cdot I.$$

Continue with the  equation $$F_s(X) = O\eqno{(2)}.$$

Here we consider $(2)$ for all $m \geq 3$. (The case $m=2$ was described in \cite{ChK}.)

Everywhere in this note, by a space we mean a metric separable space. The covering dimension of a space $X$ will be denoted by $\dim X.$ See \cite[p. 385]{E} for the definition of $\dim$ and its main properties, such as:

\begin{itemize}
\item[(i)] {\it The countable sum theorem} \cite[Theorem 7.2.1]{E}: if a  $X$ is the union  $\cup_{i=1}^\infty X_i$, where each $X_i$ is a closed subset of $X$  with $\dim X_i \leq n, n \geq 0,$  then $\dim X \leq n,$ and

\item[(ii)]{\it The monotone theorem} \cite[Theorem 7.3.4]{E}: for every subspace $A$ of a  space $X$ we have $\dim A \leq \dim X$.
\end{itemize}

Recall  that  $\dim \mathbb R^k = k$ and for any non-empty  space $X$ of cardinality less than continuum we have $\dim X = 0$.

It is easy to see that the solution set $S$ of $(2)$ (and of $(1)$ as well) is a closed subset of $M_m(\mathbb R)$ and hence $\dim S \leq m^2$. In fact, we can easily improve this estimate (similarly to \cite[Proposition 1.2]{ChK}) as follows.

\begin{prop} $\dim S \leq m^2-1$.

(The same is valid for the solution set of $(1)$.)
\end{prop}

In \cite[Theorem 3.4]{ChK} we described the solution set of $(2)$ for $m=2$ and showed (see \cite[Theorem 4.5]{ChK} ) that its dimension is
equal to $2$ for all $n \geq 2$.

Thus the following question appears naturally.

\begin{question}(\cite[Question 7.2]{ChK}) \label{question_2} Let $m \geq 3$. What is the structure of $S$? What is the dimension of $S$?
\end{question}

\section{Preliminaries I}

For a given  $B\in M_m(\mathbb R)$, define
a map
$$Conj_B\colon GL_m (\mathbb R) \to M_m(\mathbb R)$$
by $Conj_B(X) = X^{-1} B X$.

The following is easy to verify.

\begin{lem}\label{fibres} The map $Conj_B$ is continuous, and the image $Conj_B(GL_m (\mathbb R))$ is $\sigma$-compact.
\end{lem}

For a given closed subgroup $G$ of $GL_m (\mathbb R)$ denote by $GL_m (\mathbb R)/G$ the quotient space of
the right cosets of $G$ in $GL_m (\mathbb R)$, and by $\pi: GL_m (\mathbb R) \to GL_m (\mathbb R)/G$
the canonical mapping.

Since the topological group $GL_m (\mathbb R)$ is locally compact, by the use of a well known fact (see, for example, \cite[Theorem 2]{P}) we have the following
equality:
$$\dim GL_m (\mathbb R) = \dim G + \dim GL_m (\mathbb R)/G.$$
Let $G_B = \{A  \in GL_m (\mathbb R): AB = BA\}$. It is easy to see that $G_B$ is closed subgroup of
$GL_m (\mathbb R)$.

The following is evident.

\begin{prop}
\begin{itemize}
\item[(i)] Let $X, Y \in GL_m (\mathbb R)$. If $Conj_B(X) = Conj_B(Y)$ then $Y = A X$ for some $A \in G_B$.
\item[(ii)] For any $X \in GL_m (\mathbb R)$ and $A \in G_B$ we have  $Conj_B(X) = Conj_B(AX)$.
\item[(iii)] There exist a continuous bijection $\alpha_B: GL_m (\mathbb R)/G_B \to Conj_B(GL_m (\mathbb R))$
such that $Conj = \alpha_B  \pi$.
\end{itemize}
\end{prop}

Since  the space $GL_m (\mathbb R)/G_B$ is  $\sigma$-compact and normal (see \cite[Theorem 1]{P}), we get

\begin{cor}\label{dimension} $\dim Conj_B(GL_m (\mathbb R)) = dim GL_m (\mathbb R)/G_B = \dim GL_m (\mathbb R)- \dim G_B.$
\end{cor}

Recall that $(m\times m)$-matrices $B_1$ and $B_2$ with entries from a field $\mathbb F$ are {\it similar} over $\mathbb F$ if there exists an  invertible $(m\times m)$-matrix  $C$ with entries from $\mathbb F$ such that $B_1 = C^{-1}B_2C$, and the similarity is an equivalence relation.

The following is evident.
\begin{prop} For any matrices $B_1$ and $B_2$ from $M_m (\mathbb R)$
we have the sets $Conj_{B_1}(GL_m (\mathbb R))$ and $Conj_{B_2}(GL_m (\mathbb R))$ coincide
or they are disjoint.

In particular, for any $B\in M_m(\mathbb R)$ we have $Conj_{B}(GL_m (\mathbb R)) = Conj_{J}(GL_m (\mathbb R))$, where $J$ is a real Jordan canonical form of $B$.
\end{prop}


\section{Preliminaries II}


For distinct real numbers $a, b, c$ and  real numbers $p, q$, where $q \ne 0$, let 

$$J(a,b,c) = \left(\begin{array}{rrr}  a & 0 & 0 \\ 0 & b & 0  \\ 0 & 0 & c \end{array}\right),\,
J(a,a,b) = \left(\begin{array}{rrr}  a & 0 & 0 \\ 0 & a & 0  \\ 0 & 0 & b \end{array}\right),$$
$$J_1(a,a,b) = \left(\begin{array}{rrr}  a & 1 & 0 \\ 0 & a & 0  \\ 0 & 0 & b \end{array}\right),\,
J(a,a,a) = \left(\begin{array}{rrr}  a & 0 & 0 \\ 0 & a & 0  \\ 0 & 0 & a \end{array}\right),$$
$$J_1(a,a,a) = \left(\begin{array}{rrr}  a & 1 & 0 \\ 0 & a & 0  \\ 0 & 0 & a \end{array}\right),\,
J_2(a,a,a) = \left(\begin{array}{rrr}  a & 1 & 0 \\ 0 & a & 1  \\ 0 & 0 & a \end{array}\right),$$
$$J_c(a,p,q) = \left(\begin{array}{rrr}  a & 0 & 0 \\ 0 & p & q  \\ 0 & -q & p \end{array}\right),\,
U(a, p, q) = \left(\begin{array}{rrr}  a & 0 & 0 \\ 0 & p+iq & 0  \\ 0 & 0 & p-iq \end{array}\right).$$

The following is easy to verify
\begin{prop}\label{prop_commute}
\begin{itemize}
\item[(i)] $G_{J(a,b,c)} = \{A = \left(\begin{array}{rrr}  x & 0 & 0 \\ 0 & y & 0  \\ 0 & 0 & z \end{array}\right) : \det A \ne 0\}$,
\\
in particular, $\dim G_{J(a,b,c)} = 3.$
\\
\item[(ii)] $G_{J(a,a,b)} = \{A = \left(\begin{array}{rrr}  x & y & 0 \\ z & t & 0  \\ 0 & 0 & s \end{array}\right) : \det A  \ne 0\}$,
\\
in particular, $\dim G_{J(a,a,b)} = 5.$
\\
\item[(iii)] $G_{J_1(a,a,b)} = \{A = \left(\begin{array}{rrr}  x & y & 0 \\ 0 & x & 0  \\ 0 & 0 & z \end{array}\right) : \det A \ne 0\}$,
\\
in particular, $\dim G_{J_1(a,a,b)} = 3.$
\\
\item[(iv)] $G_{J(a,a,a)} =GL_3(\mathbb R)$,
in particular, $\dim G_{J(a,a,a)} = 9.$
\\
\item[(v)] $G_{J_1(a,a,a)} = \{A=\left(\begin{array}{rrr}  x & y & z \\ 0 & x & 0  \\ 0 & s & t \end{array}\right) : \det A \ne 0\}$,
\\
in particular, $\dim G_{J_1(a,a,a)} = 5.$
\\
\item[(vi)] $G_{J_2(a,a,a)} = \{A =\left(\begin{array}{rrr}  x & y & z \\ 0 & x & y  \\ 0 & 0 & x \end{array}\right) : \det A \ne 0\}$,
\\
in particular, $\dim G_{J_2(a,a,a)} = 3.$
\\
\item[(vii)] $G_{J_c(a,p,q)} = \{A=\left(\begin{array}{rrr}  x & 0 & 0 \\ 0 & y & z  \\ 0 & -z & y \end{array}\right) : \det A \ne 0\}$,
\\
in particular, $\dim G_{J_c(a,p,q)} = 3.$
\end{itemize}
\end{prop}

\section{The structure and the dimension of $S$, the case $m=3$.}

The following propositions are standard facts from matrix algebra (see, for example, \cite{L}).

\begin{prop} \label{auxprop_1} For any matrix $B \in M_3(\mathbb R)$ there exists
 a matrix $C\in GL_3(\mathbb R)$ such that $B = C^{-1} J C$, where $J$ is precisely  one of
 the following matrices:

$J(a,b,c),
J(a,a,b),
J_1(a,a,b),
J(a,a,a),
J_1(a,a,a),
J_2(a,a,a),
J_c(a,p,q).$

\end{prop}

\begin{prop}\label{auxprop_3} Any $B\in M_3(\mathbb R)$ is similar to $J_c(a,p,q)$ over $\mathbb R$ iff $B$  is similar to $U(a,p,q)$ over $\mathbb C$.
\end{prop}

Let  $f_s(x) = x^n + a_{n-1}x^{n-1} + \dots + a_1 x + a_0$ be a polynomial with real coefficients $a_i, 0 \leq i \leq n-1,$ corresponding to the polynomial  $F_s(X)$.

Consider the equation $$f_s(x) = 0 .\eqno{(3)}$$

The following lemmas are evident.
\begin{lem} \label{lem_1}
\phantom{*}

$$F_s(J(a,b,c)) = \left(\begin{array}{rrr}  f_s(a) & 0 & 0 \\ 0 & f_s(b) & 0  \\ 0 & 0 & f_s(c) \end{array}\right);\,
F_s(J(a,a,b)) = \left(\begin{array}{rrr}  f_s(a) & 0 & 0 \\ 0 & f_s(a) & 0  \\ 0 & 0 & f_s(b) \end{array}\right);$$

$$F_s(J_1(a,a,b)) = \left(\begin{array}{rrr}  f_s(a) & f'_s(a) & 0 \\ 0 & f_s(a) & 0  \\ 0 & 0 & f_s(b) \end{array}\right);\, F_s(J(a,a,a)) = \left(\begin{array}{rrr}  f_s(a) & 0 & 0 \\ 0 & f_s(a) & 0  \\ 0 & 0 & f_s(a) \end{array}\right);$$

$$F_s(J_1(a,a,a)) = \left(\begin{array}{rrr}  f_s(a) & f'_s(a) & 0 \\ 0 & f_s(a) & 0  \\ 0 & 0 & f_s(a) \end{array}\right);\,
F_s(J_2(a,a,a)) = \left(\begin{array}{rrr}  f_s(a) & f'_s(a) & \frac{1}{2}f''_s(a) \\ 0 & f_s(a) & f'_s(a)  \\ 0 & 0 & f_s(a) \end{array}\right);$$

$$F_s(U(a,p,q)) = \left(\begin{array}{rrr}  f_s(a) & 0 & 0 \\ 0 & f_s(p+iq) & 0  \\ 0 & 0 & f_s(p-iq) \end{array}\right).$$

In particular,

$a, b, c$ are distinct real roots of $(3)$ iff $J(a,b,c)$ is a solution of $(2)$;

$a, b$ are distinct real roots of $(3)$  iff $J(a,a,b)$ is a solution of $(2)$;

$a, b$ are distinct real roots of $(3)$ and $a$ is a root of multiplicity $\geq 2$ iff
$J_1(a,a,b)$ is a solution of $(2)$;

$a$ is a real root of $(3)$ iff $J(a,a,a)$ is a solution of $(2)$;

$a$ is a real root of $(3)$ of multiplicity $\geq 2$ iff $J_1(a,a,a)$ is a solution of $(2)$;

$a$ is a real root of $(3)$ of multiplicity $\geq 3$ iff $J_2(a,a,a)$ is a solution of $(2)$;

$a$ is a real root of $(3)$ and $p \pm iq$ are non-real roots of $(3)$ iff
$F_s(U(a,p,q)) = O$.
\end{lem}

\begin{lem} \label{lem_2} If $(m\times m)$-matrices $A$ and $B$ (with real or complex entries) are similar over $\mathbb R$ or $\mathbb C$, and $F_s(A) = O,$ then $F_s(B) = O.$

In particular, $J_c(a,p,q)$ is a solution of $(2)$ iff $a$ is a real root of $(3)$ and $p\pm iq$ are non-real roots of $(3)$.
\end{lem}

Let $S(B) =Conj_{B}(GL_3(\mathbb R))$ for any $B \in M_3(\mathbb R)$.

Using Propositions \ref{auxprop_1}, \ref{auxprop_3}, and Lemmas \ref{lem_1}, \ref{lem_2} we get the following.

\begin{thm} \label{character_theorem}The solution set $S$ of the equation $(2)$ is either empty or a disjoint union of the following sets:

$S(J(a,b,c))$, whenever $a,b,c$ are distinct real roots of  $(3)$,

$S(J(a,a,b))$, whenever $a,b$ are distinct real roots of  $(3)$,

$S(J_1(a,a,b))$, whenever $a, b$ are distinct real roots of  $(3)$ and
$a$ is a root of multiplicity $\geq 2$,

$S(J(a,a,a))$, whenever $a$ is a real root of  $(3)$,

$S(J_1(a,a,a))$, whenever $a$ is a real root of $(3)$ of multiplicity $\geq 2$,

$S(J_2(a,a,a))$, whenever $a$ is a real root of $(3)$ of multiplicity $\geq 3$,

$S(J_c(a,p,q))$, $q > 0$, whenever $a$ is a real root of $(3)$ and
$p \pm iq$ are non-real roots of $(3)$.
  \end{thm}

Using Proposition \ref{prop_commute} and Corollary \ref{dimension} we get the following corollary.
\begin{cor}
\begin{itemize}
\item[(i)] $\dim S(J(a,b,c)) = \dim S(J_1(a,a,b)) =$

$\dim S(J_2(a,a,a)) = \dim S(J_c(a,p,q)) = 6.$

\item[(ii)] $\dim S(J(a,a,b)) = \dim S(J_1(a,a,a)) = 4$ .

\item[(iii)]  $\dim S(J(a,a,a))= 0$,
moreover,

$S(J(a,a,a)) = \{\left(\begin{array}{rrr}  a & 0 & 0 \\ 0 & a & 0  \\ 0 & 0 & a \end{array}\right)\}$.
\end{itemize}
\end{cor}

 \begin{cor}\label{cor_sigma_compact} The solution set $S$ of the equation $(2)$, $m=3$, is
 $\sigma$-compact, and $\dim S$ is equal to $-1$ or it is the maximum of dimensions of sets $S(J)$ from Theorem \ref{character_theorem}, i.e. $0, 4$ or $6$.
  \end{cor}


\section{Some examples of solution sets, the case $m=3$.}

Here we will describe the solution sets of some equations. All facts here follow from the results of Section~4.

\begin{example} The solution set $S_1$ of $X^2+ I =O$ is  empty, and so
$\dim S_1 = -1.$

\end{example}

\begin{example}
\begin{itemize}
\item [(i)] The solution set $S_2$ of $X =O$ is equal to $S(J(0,0,0)) = \{O\}$, and so $\dim S_2 = 0.$
\item [(ii)] The solution set $S_3$ of $X - I=O$ is equal to $S(J(1,1,1)) = \{I\}$, and so $\dim S_3 = 0.$
\end{itemize}
\end{example}

\begin{example} The solution set $S_4$ of $X^2 =O$ is equal to $S(J(0,0,0)) \cup S(J_1(0,0,0))$, and
so $\dim S_4 = 4.$
\end{example}

\begin{example} The solution set $S_5$ of $X^2(X-I) =O$ is equal to $S(J(0,0,0)) \cup S(J_1(0,0,0)) \cup
S(J(1,1,1)) \cup S(J(0,0,1)) \cup S(J(1,1,0)) \cup S(J_1(0,0,1))$,
and
so $\dim S_5 = 6.$

Moreover, $S_5 \setminus (S_3 \cup S_4) = S(J(0,0,1)) \cup S(J(1,1,0)) \cup S(J_1(0,0,1))$,

 and so
$\dim S_5 \setminus (S_3 \cup S_4) = 6 > 4 = \dim (S_3 \cup S_4).$
\end{example}

\begin{example} The solution set $S_6$ of $X^2 - I=O$ is equal to
$S(J(1,1,1)) \cup S(J(-1,-1,-1)) \cup S(J(1,1,-1)) \cup S(J(-1,-1,1))$, and so $\dim S_6 = 4.$
\end{example}

\begin{example} The solution set $S_7$ of $X(X^2-I) =O$ is equal to $S(J(0,0,0)) \cup
S(J(1,1,1)) \cup  S(J(-1,-1,-1)) \cup S(J(0,0,1)) \cup S(J(1,1,0)) \cup S(J(0,0,-1)) \cup S(J(-1,-1,0))
\cup S(J(1,1,-1)) \cup S(J(-1,-1,1)) \cup S(J(-1,0,1))$, and
so $\dim S_7 = 6.$

Moreover, $S_7 \setminus (S_2 \cup S_6) = S(J(0,0,1)) \cup S(J(1,1,0)) \cup S(J(-1,0,1)) \cup S(J(0,0,-1)) \cup S(J(-1,-1,0))$,

 and so
$\dim S_7 \setminus (S_2 \cup S_6) = 6 > 4 = \dim (S_2 \cup S_6).$
\end{example}

\begin{example} The solution set $S_8$ of $X(X^2+I) =O$ is equal to $S(J(0,0,0)) \cup
S(J_c(0,0,1)) $, and
so $\dim S_8 = 6.$

Moreover, $S_8 \setminus (S_1 \cup S_2) =S(J_c(0,0,1)) $,

 and so
$\dim S_8 \setminus (S_1 \cup S_2) = 6 > 0 = \dim (S_1 \cup S_2)$
\end{example}

\begin{example} The solution set $S_9$ of $X^3 =O$ is equal to $S(J(0,0,0)) \cup S(J_1(0,0,0)) \cup
S(J_2(0,0,0))$, and
so $\dim S_9 = 6.$

Moreover, $S_9 \setminus (S_2 \cup S_4) =S(J_2(0,0,0)) $,

 and so
$\dim S_9 \setminus (S_2 \cup S_4) = 6 > 4 = \dim (S_2 \cup S_4).$
\end{example}

\section{The structure of $S$ for $m \geq 4$}
Let the equation $(3)$ have  distinct real roots $r_1, \dots, r_k$ and distinct non-real roots
$c_1, \dots, c_l$. Using these two sets of numbers one can construct only finitely many non-similar real Jordan canonical forms
$J_1, \dots, J_p$ which are elements of $M_m(\mathbb R)$.

Let $S$ be  the solution set of the corresponding equation $(2)$. Then $S$ is the disjoint union 
$\cup_{i=1}^p Conj_{J_i}(GL_m(\mathbb R))$ which is $\sigma$-compact, and $\dim S$ is $-1$ or
 the maximum of $\dim Conj_{J_i}(GL_m(\mathbb R)), i \leq p$.

\begin{prop} If $(3)$ has only distinct real roots and the degree of $(3)$ is greater or equal to $m$ then the dimension of $S$ is $m^2-m$.
\end{prop}

\begin{prop} If $(3)$ has only non-real roots and $m = 2k+1, k \geq 1,$ then the dimension of $S$ is $-1$.
\end{prop}

\section{Concluding remarks}

\begin{question} Let $m \geq 4$. What values can $\dim S$ admit?
\end{question}

\vskip 0.5 cm
\noindent(V.A. Chatyrko)\\
Department of Mathematics, Linkoping University, 581 83 Linkoping, Sweden.\\
vitalij.tjatyrko@liu.se

\vskip0.3cm
\noindent(A. Karassev) \\
Department of Mathematics,
Nipissing University, North Bay, Canada\\
alexandk@nipissingu.ca


\begin{thebibliography}{99}
\bibitem[ChK]{ChK} V.~A.~Chatyrko, A.~Karassev, Simple polynomial equations over $(2 \times 2)$-matrices,
arXiv:2506.07689

\bibitem[E]{E} R.~Engelking, General Topology, Heldermann Verlag, Berlin, 1989.

\bibitem[FS]{FS} D.~Fuchs, A.~Schwarz,  A matrix Vieta Theorem, E. B. Dynkin Seminar, Amer. Math. Soc. Thansl. Ser. 2, 169 (1996)

\bibitem[G]{G} S.~Gelfand, On the number of solutions of a quadratic equation, In: Globus: General Mathematical Seminar, 1,
Independent University of Moscow, Moscow, 2004, 124-133 (in Russian)

\bibitem[P]{P} B.~A.~Pasynkov, On the coincidence of various definitions of dimensionality for factor spaces of locally bicompact groups, Uspekhi Matematicheskikh Nauk, 1962, 17(5),  129–135	

\bibitem[L]{L} H.~Lutkepohl, Handbook of Matrices, John Wiley \& Sons, 1996

\bibitem[S]{S} M.~Slusky, Zeros of $(2 \times 2)$-matrix polynomials,  Comm. Algebra, 2010,  38(11),  4212-4223

\bibitem[W]{W} R.~L.~Wilson, {\em Polynomial equations over matrices}, Rutgers University, manuscript

\end{thebibliography}
\end{document}